\documentclass{article}

\usepackage[english]{babel}

\usepackage[letterpaper,top=2cm,bottom=2cm,left=3cm,right=3cm,marginparwidth=1.75cm]{geometry}

\usepackage{amsfonts}
\usepackage{amsmath}
\usepackage{amsthm}
\usepackage{amssymb}
\usepackage{graphicx}
\usepackage[colorlinks=true, allcolors=blue]{hyperref}
\usepackage{float}
\usepackage{booktabs}
\usepackage[table,dvipsnames]{xcolor}

\usepackage{stmaryrd}
\usepackage{subcaption}
\usepackage{comment}
\usepackage{authblk}

\newcommand{\R}{\mathbb{R}}
\newcommand{\N}{\mathbb{N}}
\renewcommand{\P}{\mathbb{P}}
\newcommand{\interv}[1]{\llbracket #1\rrbracket}
\newcommand{\vc}[1]{{\bf{#1}}}

\newcommand{\reviewerA}[1]{#1}

\date{}

\title{Hermite Semi-Lagrangian schemes on triangular meshes for advection equations}
\author[1,2,3]{Ali Elarif}
\author[4]{Michel Mehrenberger}
\author[1,2]{Laurent Navoret}

\affil[1]{IRMA, UMR 7501 Université de Strasbourg et CNRS, 7 rue René Descartes, 67000 Strasbourg, France}
\affil[2]{Inria, IRMA, Université de Strasbourg, CNRS UMR 7501, 7 rue René Descartes, 67084 Strasbourg, France}
\affil[3]{Nodji Consulting, 236B rue du Faubourg des postes, 59000 Lille}
\affil[4]{Aix Marseille Université, CNRS, I2M, Marseille, France}

\begin{document}

\maketitle

\begin{abstract}
High-order Hermite semi-Lagrangian schemes on unstructured triangular grids are proposed for advection equations, based on Bell and Argyris finite elements. Nodal semi-Lagrangian schemes transport point values together with gradients and Hessians, while projection semi-Lagrangian schemes are basically semi-Lagrangian Discontinuous Galerkin schemes (SLDG) with straight remaining backward triangles and shared degrees of freedoms. Stability and convergence studies of the different schemes are carried out numerically.
\end{abstract}

\section{Introduction}
\label{sec:introduction}

Semi-Lagrangian (SL) schemes are a class of schemes widely used to solve advection equations. At each time iteration, the approximate solution is first advected exactly, i.e. expressed as the solution at the previous time step evaluated at the foot of the characteristics, and then \reviewerA{pushed} back to the approximation space by nodal evaluations (nodal SL) or by $L^2$ projection (projection SL). These schemes have the advantage of not being constrained by a CFL-type condition, which is of great advantage in the presence of unnecessarily  small mesh cells. However, while high-order semi-Lagrangian schemes have been proposed on uniform Cartesian (e.g. odd degree Lagrange nodal SL \cite{charles2013enhanced,ferretti2020stability,yang2021highly}) or non-uniform Cartesian grids (Hermite nodal SL, Semi-Lagrangian Discontinuous Galerkin), the extension to unstructured meshes is more challenging due to stability and computational issues.

Regarding nodal SL schemes on unstructured triangular grids, one popular solution is to use $\P_2$  interpolation \cite{ferretti2007stability}, but $\P_3$ or higher order interpolation are unstable 
\cite{crouseilles2011discontinuous}.
Another approach is to consider a larger stencil together with least square reconstruction \cite{carre2009polynomial,bergami2020high}. 
Finally, the third option consists in using Hermite interpolation by transporting nodal values of the function and its gradients \cite{besse2003semi} and based on the reduced Hsieh-Clough-Tucker (rHCT) approximation space. Note that in \cite{besse2003semi}, only constant advection fields have been considered.

Regarding projection SL schemes, also called Semi-Lagrangian Discontinuous Galerkin schemes (SLDG) when combined with discontinuous polynomial approximation, stability issues appear to be ruled out, as the method is stable due to the $L^2$ projection. However, the method requires the computation of integrals involving quantities defined on backward cells and their approximations can lead to instabilities \cite{restelli2006semi,pironneau1982transport}.
In order to prevent from such instabilities, further projection semi-Lagrangian schemes are based on mesh intersection algorithms and has been developed on Cartesian or polar grids 
\cite{lauritzen2010conservative,cai2021high,crouseilles2014new,cai2019high,einkemmer2024semi}.
In the context of triangular meshes, a recent work proposes a SLDG scheme combined with a Runge-Kutta Discontinuous Galerkin method \cite{cai2022eulerian}, where backward cells do not need to be curved to be still high order accurate. 
We also refer to \cite{farrell2011conservative,alauzet2010p1} for remapping schemes on triangular meshes, based on mesh intersections.


In this work, we propose to explore the use of high-order Hermite elements on triangles, i.e. Bell and Argyris elements \cite{bell1969refined,okabe1994explicit}. In the context of nodal SL schemes,  this requires to also advect the nodal values of the Hessian. Moreover, the projection SL method can also be adapted to other discretizations than Discontinuous Galerkin: we will also explore the use of the Hermite elements in this context, which permit in particular to share some degrees of freedom. Note that, at the price of some loss of precision, the backward triangles will be approximated by straight and not curved triangles \cite{carre2009polynomial}. The outline of this paper is as follows. In Section 2, we describe the schemes and numerical results are given in Section 3.

\section{Hermite Semi-Lagrangian schemes}
\label{sec:1}

We are interested in solving the following two-dimensional advection problem:
\begin{align}
    &\partial_t \rho + a \cdot \nabla_x \rho = 0,\\
    &\rho(0,x) = \rho_{\text{init}}(x),
\end{align}
where $\rho(t,x) \in \R$, with $t \geq 0$ and $x \in \Omega \subset \R^2$, denotes the advected density, $a(t,x) \in \R^2$ the advection velocity field and $\rho_{\text{init}}(x) \in \R$ the initial density. The equation is supplemented either with Dirichlet conditions at the inflow boundary. The solution to this equation is constant along the characteristic curves. In particular, given discrete times $t^n = n \Delta t$, with $\Delta t \geq 0$, we have:
\begin{equation*}
    \rho(t^{n+1},x) = \rho\left(t^n, X(t^n; t^{n+1}, x)\right)
\end{equation*}
where $X(s;t,x) \in \R^2$ denotes the characteristic curves associated with the advection equation, solutions to the following set of differential equations:
\begin{align*}
    &\partial_s X(s;t,x) = a(s, X(s;t,x)),\\
    &X(t;t,x) = x.
\end{align*} 
Given a triangular mesh $\mathcal{T}_h$ of the domain $\Omega$, with maximal cell diameter $h$, we aim at constructing an approximate solution $\rho_h^n$, at each discrete time $t^n = n \Delta t$, with a given time step $\Delta t > 0$. We will use Hermite finite element space approximation.

\subsection{Approximate solution: Hermite finite element space}

On each triangle element $K$ of the mesh, we consider the following degrees of freedom $\Sigma$ and approximation space $P$.
\begin{description}
    \item[\textbf{reduced Hsieh-Clough-Tocher (rHCT)}] The degrees of freedom are the nodal values and the two partial derivatives at each node:
    \[\Sigma = \left\{\rho(x_i),\ \partial_{x_k}\rho(x_i),\quad \forall i\in\interv{1,3},\forall k \in\interv{1, 2}\right\}.\]
    Using these $9$ degrees of freedom, the approximation space is defined by:
    \[
    P = \left\{ p \in \mathcal{C}^1(K),\quad p_{|K_i} \in \P_3,\quad \partial_n p_{|\partial K_i\cap\partial K} \in \P_1,\quad \forall i \in \interv{1,3} \right\},
    \]
    where $(K_i)$ are the three subtriangles whose vertices are the barycenter and the vertices of $K$. This space contains the polynomials of degree $2$ (space of dimension $6$) and the associated interpolation error is thus locally of order $3$. 
\end{description}
The global approximation space is thus composed of piecewise polynomials of degree 3 (on subtriangles) with $\mathcal{C}^1$ global regularity. 
We are also interested in higher-order Hermite finite element approximation.
\begin{description}
    \item[\textbf{Hsieh-Clough-Tocher (HCT)}] Compared with rHCT, $3$ additional degrees of freedom are considered. These are the normal derivatives at the middle of the edges: \[\Sigma = \{\rho(x_i),\ \partial_{x_k}\rho(x_i),\ \partial_{n}\rho((x_i+x_j)/2),\ \forall i \neq j \in \interv{1,3}, \forall k\in \interv{1, 2}\}.\]   Using these $12$ degrees of freedom, the approximation space is now defined by:
    \[
    P = \left\{ p \in \mathcal{C}^1(K),\quad p_{|K_i} \in \P_3,\quad \forall i \in \interv{1,3} \right\},
    \]
    where $(K_i)$ are still the three subtriangles whose vertices are the barycenter and the vertices of $K$.
    Polynomials of degree $3$ (space of dimension 
    $10$) are exactly reproduced and the associated interpolation error is thus locally of order $4$.
    \item[\textbf{Bell}] Instead of considering derivatives at the middle of the edges, Bell element uses instead the three second order derivatives at each three nodes. The degrees of freedom are thus given by:
    \[\Sigma = \{\rho(x_i),\ \partial_{x_k}\rho(x_i),\ \partial_{x_k x_\ell}\rho(x_i),\ \forall i \in \interv{1,3},\forall k,\ell\times\interv{1, 2}^2\}.\] 
    Using these $18$ degrees of freedom, the associated approximation space is given by:
    \[
    P = \left\{ p \in \P_5,\quad \partial_n p_{|\partial K} \in \P_3 \right\}.
    \]
    Polynomials of degree 4 (space of dimension 15) are exactly reproduced and the associated interpolation error is thus locally of order 5.
    \item[\textbf{Argyris}] Finally, the Argyris element considers three additional degrees of freedom, which are the normal derivatives at the middle of the edges:
    \begin{align*}
    \Sigma = \{\rho(x_i),\ \partial_{x_k}\rho(x_i),\ & \partial_{x_k x_\ell}\rho(x_i),\\
    &\partial_{n}\rho((x_i+x_j)/2),\ \forall i \neq j \in \interv{1,3},\forall k,\ell\times\interv{1, 2}^2\}.
    \end{align*}
    Using these $21$ degrees of freedom, the approximation space equals $P = \P_5$ and the associated interpolation error is thus locally of order $6$. 
\end{description}
We will also explore the combination of the Argyris finite element space for $\rho$ and the reduced HCT for its derivatives $\partial_{x_1}\rho$, $\partial_{x_2}\rho$. The scheme will be refered as \textbf{Argyris - grad rHCT}. The local error on the derivatives is thus of order $3$ instead of $5$: the expected overall error is thus at most $4$,  i.e., an additional order of convergence compared to rHCT.
Given the associated set of degrees of freedom, each finite element provides an interpolation operator over the whole domain $\Omega$:
\begin{equation*}
    \rho_h^n(.) = \mathcal{I}\Big(\rho_h^n(x_i),\ \partial_{x_k}\rho_h^n(x_i),\ \partial_{x_k x_\ell}\rho_h^n(x_i),\  \partial_{n}\rho_h^n((x_i+x_j)/2), \forall i \neq j\Big)
\end{equation*}
By construction, this approximate solution is piecewise polynomial (on subtriangles for HCT elements) with $C^1$ global regularity. Table~\ref{tab:hermitefem} provides a summary of these elements.

\begin{table}[h]
\caption{High order Hermite finite elements}
\label{tab:hermitefem}       
%
%
\centering
\begin{tabular}{llcl}
    \toprule
    &$\Sigma$&dimension&$P$\\
    \hline\noalign{\smallskip}
    rHCT & $\rho(x_i)$,$\partial_{x_j}\rho(x_i)$ &9 & $\P_2 \subset  P \subset \P_3$\\
    HCT&$\rho(x_i)$, $\partial_{x_j}\rho(x_i)$, $\partial_{n}\rho((x_i+x_j)/2)$&12 & $\P_3 \subset  P \subset \P_4$\\
    Bell&$\rho(x_i)$,$\partial_{x_j}\rho(x_i)$,$\partial_{x_j x_k}\rho(x_i)$ &18 & $\P_4 \subset  P \subset \P_5$\\
    Argyris&$\rho(x_i)$,$\partial_{x_j}\rho(x_i)$,$\partial_{x_j x_k}\rho(x_i)$,  $\partial_{n}\rho((x_i+x_j)/2)$ &21 & $P = \P_5$\\
    \noalign{\smallskip}\bottomrule
\end{tabular}
\end{table}

\subsection{Update the degrees of freedom}

Consequently, to update the approximate solution in time, we need to update the degrees of freedom. This can be done in two ways:
\begin{description}
    \item[\textbf{(nodal SL)}]  The nodal values are updated using the following formulas:
    \begin{align*}
        \rho^{n+1}_h(x_i) &= \rho_h^n(X^n(x_i)),&&\\ 
        \partial_{x_j}\rho^{n+1}_h(x_i) &= \nabla_x \rho^n_h(X(x_i)) \cdot \partial_{x_k} X^n(x_i),&&\forall k \in \{1,2\},\\ 
        \partial_{x_k x_\ell}\rho^{n+1}_h(x_i) &= \nabla_x \rho^n_h(X(x_i)) \cdot \partial_{x_k x_\ell} X^n(x_i),&&\forall k,\ell \in \{1,2\},\\
        &\quad + \partial_{x_k} X^n(x_i) \cdot \left(\nabla^2 \rho_h^{n}(X^n(x_i))\, \partial_{x_\ell} X^n(x_i)\right),&&\\
        \partial_{n}\rho_h^n((x_i+x_j)/2) &= \nabla_x \rho^n_h(X^n((x_i+x_j)/2)) \cdot \partial_{n} X^n((x_i+x_j)/2)
    \end{align*}
    where $X^n(x)$ stands for $X(t^n; t^{n+1}, x)$.
    \item[\textbf{(projection SL)}] The updated density is defined by the $L^2$ projection on the approximate finite element space of the exact solution:
    \begin{equation*}
        \int_{\Omega} \rho_h^{n+1}(x)\, \varphi_i(x)\, dx =  \int_{\Omega} \rho_h^{n}(X^n(x))\, \varphi_i(x)\, dx
    \end{equation*}
   for all $\varphi_i(x)$ basis function associated with the degrees of freedom, where $X^n(x)$ still stands for $X(t^n; t^{n+1}, x)$.
\end{description}
The projection SL method directly ensures $L^2$ stability. However, \reviewerA{the stability is no more guaranteed if non-exact projection is used \cite{morton1988stability}} and the numerical scheme can also be difficult to implement in practice because the calculation of integrals requires mesh intersections. In the appendix, we propose an implementation using Firedrake \cite{farrell2011conservative, FiredrakeUserManual}. The nodal semi-Lagrangian scheme does not require such computationally expensive task, but the derivatives of the characteristic field are required.

\subsection{Approximation of the characteristic curves}

The characteristic curves are here approximated using an explicit Runge-Kutta (RK) numerical solver, whose general expression is given by:
\begin{align*}
    &\vc{X}(t^n;t^{n+1}, \vc{x}) = \vc{x} - \Delta t \sum_{i=1}^p b_i \vc{k}_i(\vc{x}),\\
    &\vc{k}_i(\vc{x}) = \vc{a}\Big(t^n - c_i \Delta t, \vc{x} - \Delta t \sum_{j=1}^{i-1} a_{i,j} \vc{k}_j(\vc{x})\Big),\quad  \forall i \in \interv{1, p},
\end{align*}
where $p \in \N$ denotes the number of steps and $(b_i)_{1\leq i \leq p}$, $(c_i)_{1\leq i \leq p}$, and $(a_{i,j})_{1\leq j < i \leq p}$ the coefficients. In practice, we will use the fourth-order RK scheme. To increase the precision of the resolution of the characteristics, $m$ sub-steps can be used with a time step $\delta t > 0$ such that $m\delta t = \Delta$. As required for the nodal SL scheme, the derivatives of the characteristic curves are detailed in the  appendix.

\section{Numerical results}

In this section, we perform a detailed study of the accuracy of the different scheme on three different test cases. 

\subsection{Translation}

\begin{figure}
    \centering
    \includegraphics[width=0.45\textwidth]{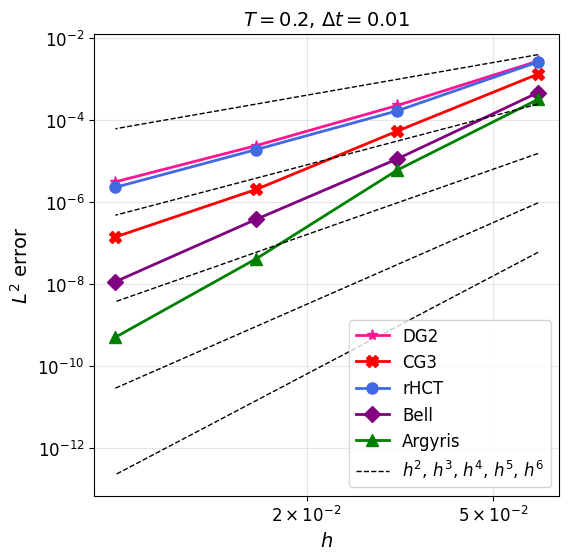}
    \includegraphics[width=0.45\textwidth]{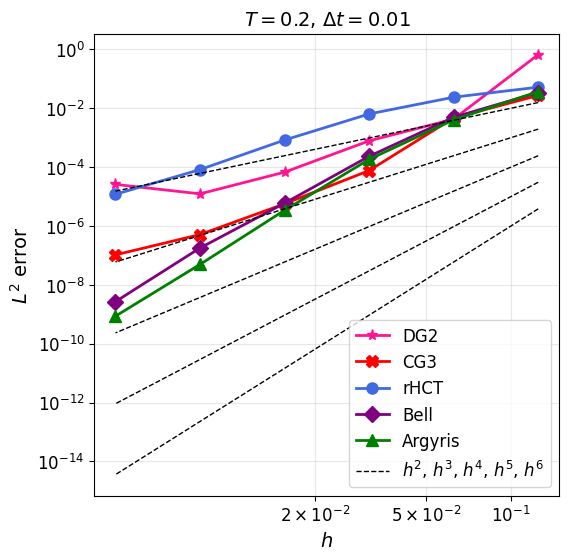}
    \caption{(Translation) Convergence curves with $\Delta t = 0.01$. Left: projection  SL. Right: nodal SL.} 
    \label{fig:translation_convergence}
\end{figure}

We consider the case of the advection with a constant advection field $a(t,x) = (1,1)^T$ in the square domain $\Omega = [0,1]^2$. The initial data is a Gaussian function:
\begin{equation}
    \rho_{\text{init}}(x) = \exp\left( - \frac{|| x - x_c||^2}{2\sigma^2}\right), 
    \label{eq:gaussian}
\end{equation}
centered at $x_c = (0.4, 0.4)^T$ and with standard deviation equal to $\sigma = 0.04$. As the characteristic are known exactly, there is no time error. The errors result only from the succession of the spatial interpolations. Figure~\ref{fig:translation_convergence} shows the convergence curves for the different nodal and projection Hermite SL schemes. We also compare these scheme with the DG2 and CG3 SL schemes, which have the same order of number of degrees of freedom (dofs) than Argyris: for a mesh with $N$ vertices and thus approximately $3N$ edges and $2N$ triangles, Argyris have $9N$ dofs while DG2 (resp. CG3) have $12N$ dofs (resp. $9N$ dofs).

\subsection{Rotation}

\begin{figure}[h]
    \centering
    \includegraphics[width=0.35\textwidth]{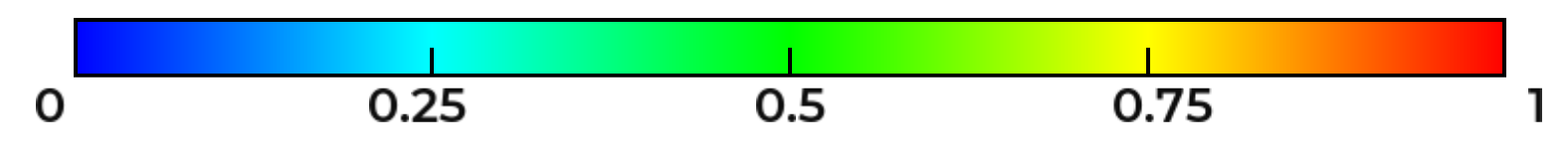}\newline
    \includegraphics[width=3.8cm]{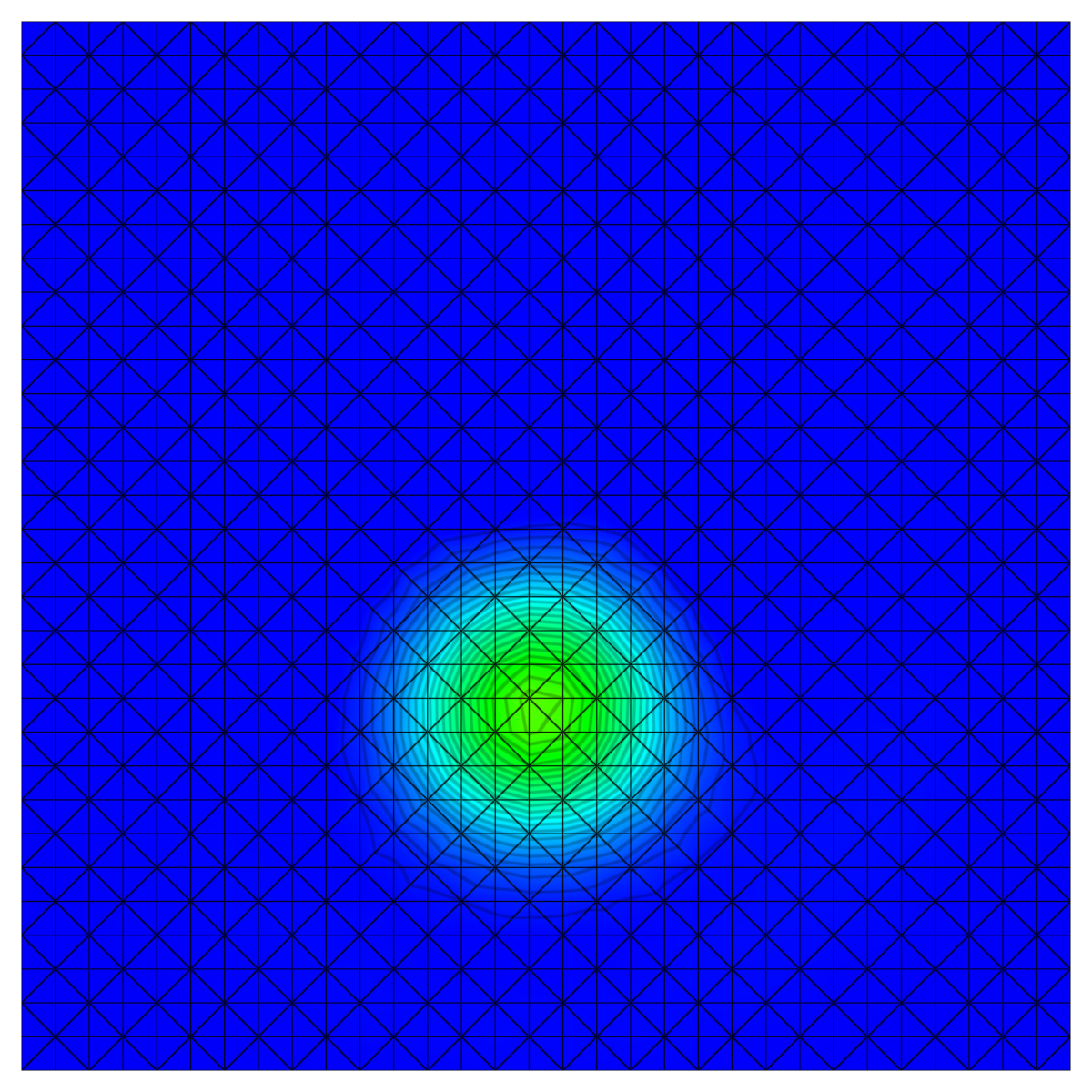}
    \includegraphics[width=3.8cm]{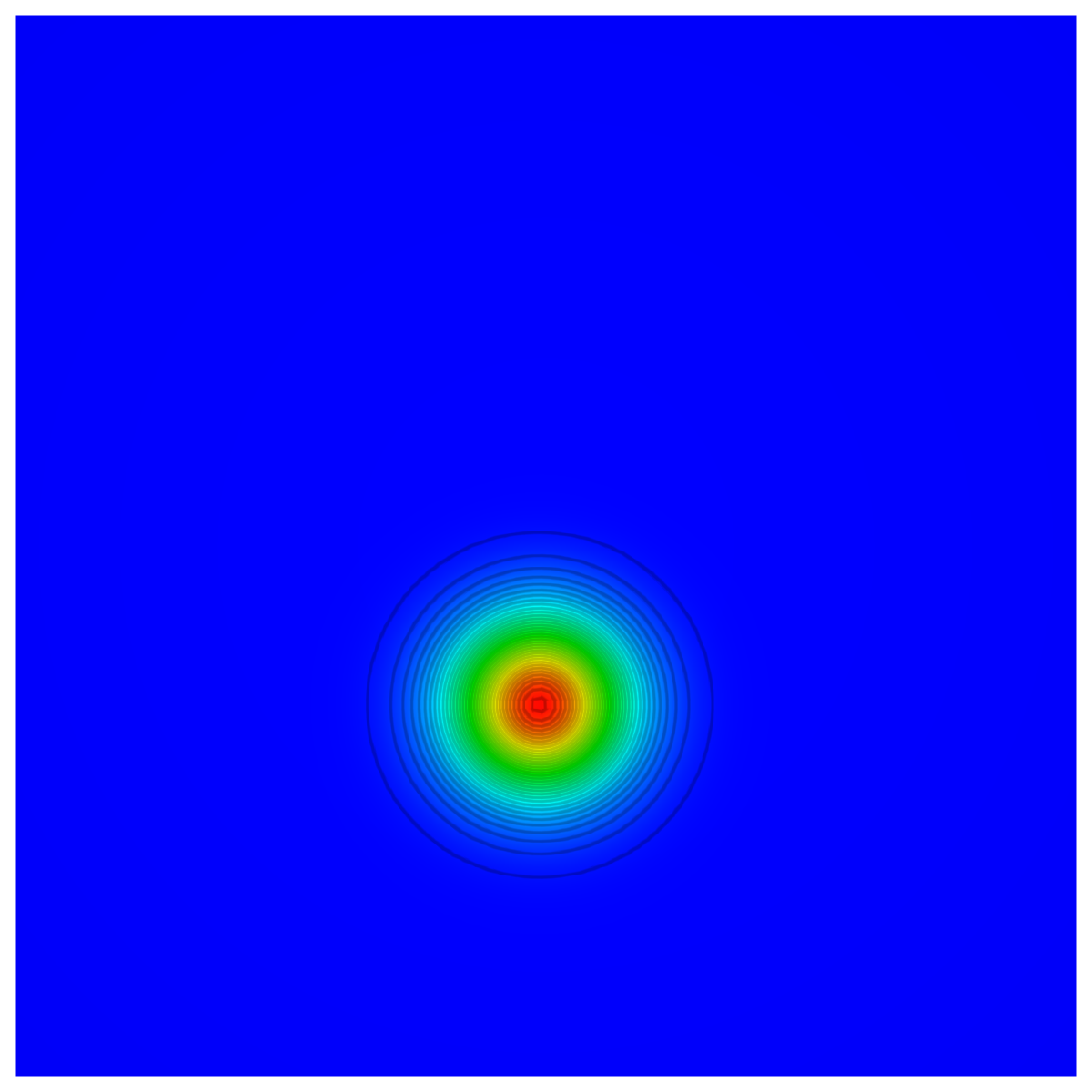}\newline
    $N=32$, $h\approx 2\times 10^{-1}$\qquad\qquad $N=128$, $h\approx 5\times 10^{-2}$
    \caption{(Rotation) Density at times $T=256$, obtained with the Argyris-gradrHCT nodal SL scheme with $\Delta t=1/4$ and two different spatial resolutions.}
    \label{fig:rotation_rho} 
\end{figure}

We consider the rotation test case, where the domain is the unit square $\Omega = [-\pi, \pi]^2$ and the advection field is given by:
\begin{equation*}
    a(t, x) = \begin{bmatrix}
        - x_2\\
        x_1
    \end{bmatrix}.
\end{equation*}
The initial condition is still the Gaussian function \eqref{eq:gaussian} but
centered at $x_c = (0.3\pi, 0)^T$ and with standard deviation  $\sigma = 0.35$. The characteristic curves are rotations:
\begin{equation*}
    X(t^{n}; t^{n+1}, x) = \begin{bmatrix}
        x_1 \cos(-\Delta t) - x_2 \sin(-\Delta t) \\
        x_1 \sin(-\Delta t) + x_2 \cos(-\Delta t) 
    \end{bmatrix}.
\end{equation*}
Like in the translation case, the errors result only from the spatial interpolations. Figure~\ref{fig:rotation_rho} 

\begin{figure}
    \centering
    \includegraphics[width=0.45\textwidth]{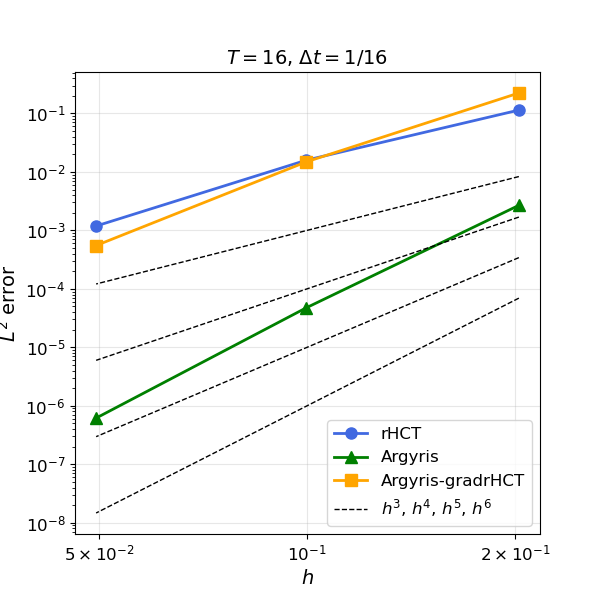}
    \includegraphics[width=0.45\textwidth]{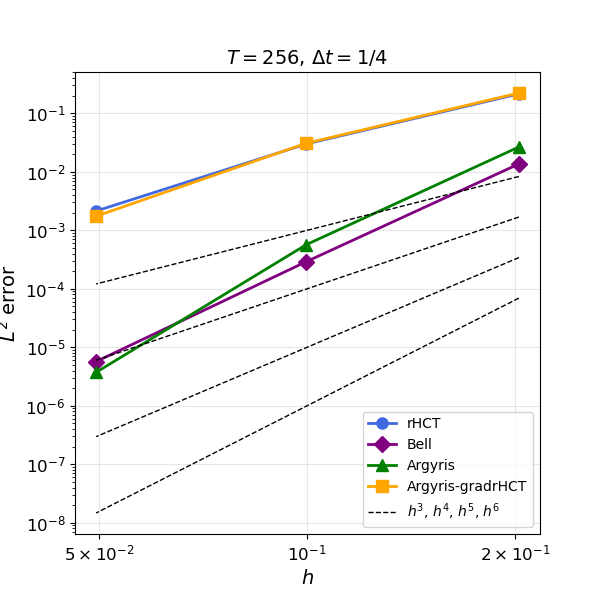}
    \caption{(Rotation) Convergence curves with two different time steps. Left: $\Delta t = 1/16$. Right: $\Delta t = 1/4$.} 
    \label{fig:rotation_convergence}
\end{figure}

\begin{figure}
    \centering
    \includegraphics[width=0.45\textwidth]{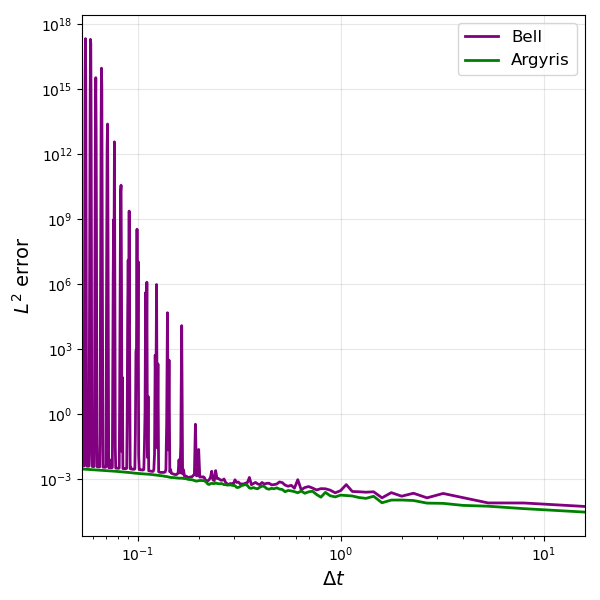}
    \caption{(Rotation test case) $L^2$ error at final time $T= 16$ as a function of \reviewerA{time step} for nodal SL schemes with $N=32$ ($h \approx 0.2$).} 
    \label{fig:stability_rotation}
\end{figure}

Figure~\ref{fig:rotation_convergence} shows the convergence curves for the Hermite nodal SL schemes for two different time steps. As expected, the Argyris scheme has order of convergence larger than $6$ and the rHCT scheme order larger than $3$. Argyris gradrHCT has a slightly better accuracy than rHCT. Although the number of dofs of Argyris \reviewerA{is} $3$ times larger than for rHCT ($9N$ instead of $3N$), Argyris is more competitive, due to the higher order of accuracy. We note that Bell scheme does not work when the time step equals $1/4$ (Fig.\ref{fig:rotation_convergence}). Figure~\ref{fig:stability_rotation} actually shows that, unlike Argyris, Bell develop instabilities for some values of time steps.  

\subsection{Swirling deformation flow}

\begin{figure}
    \centering
    \includegraphics[width=0.35\textwidth]{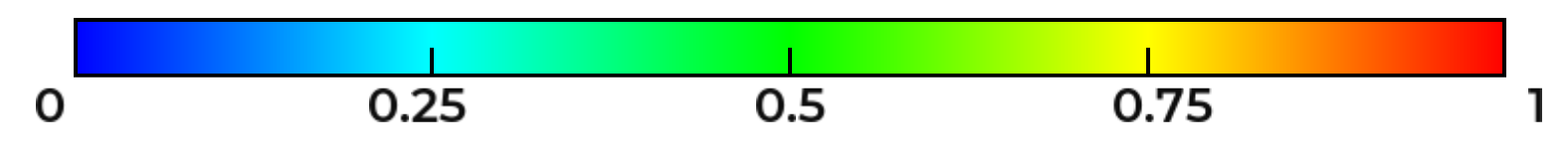}\newline
    \includegraphics[width=0.35\textwidth]{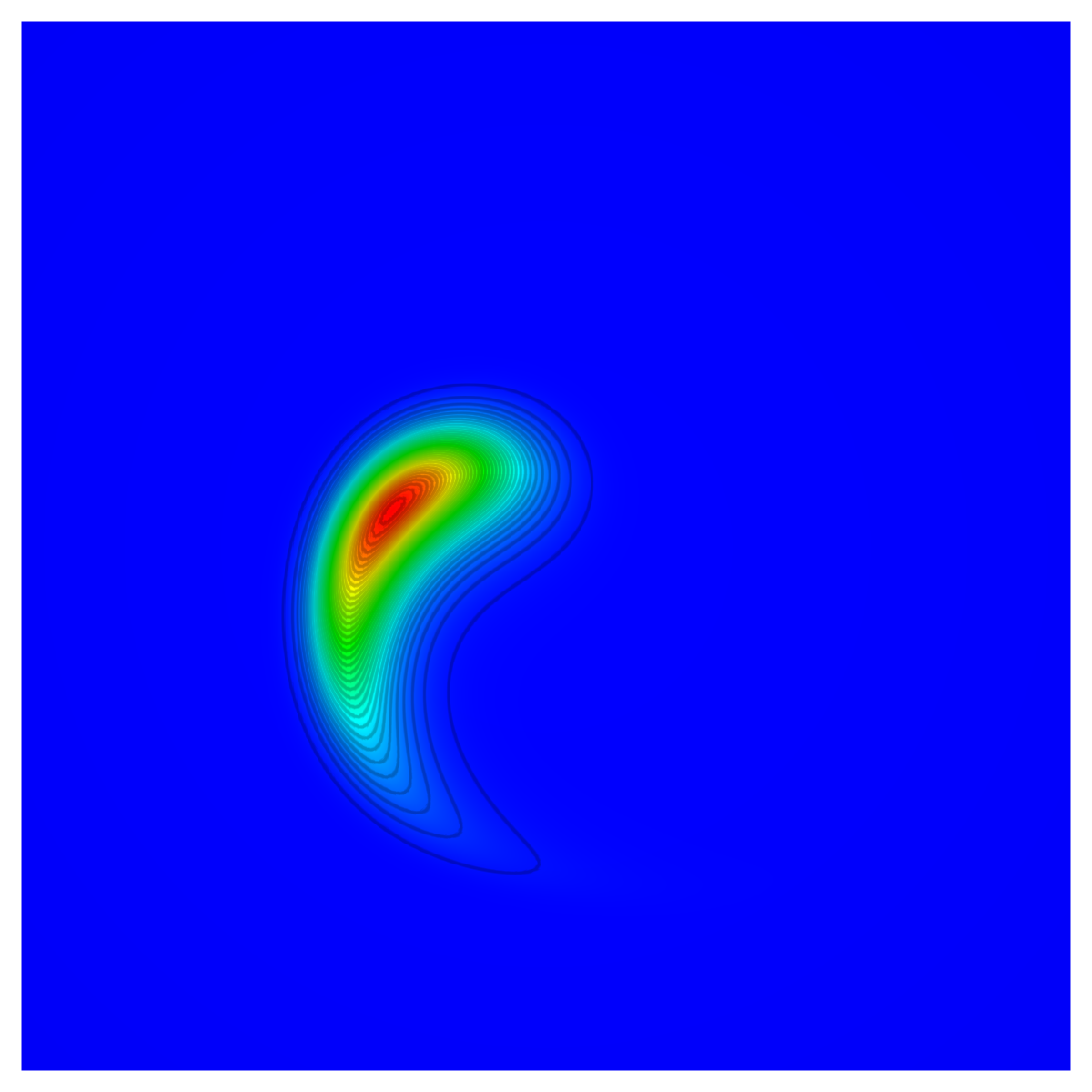}
    \includegraphics[width=0.35\textwidth]{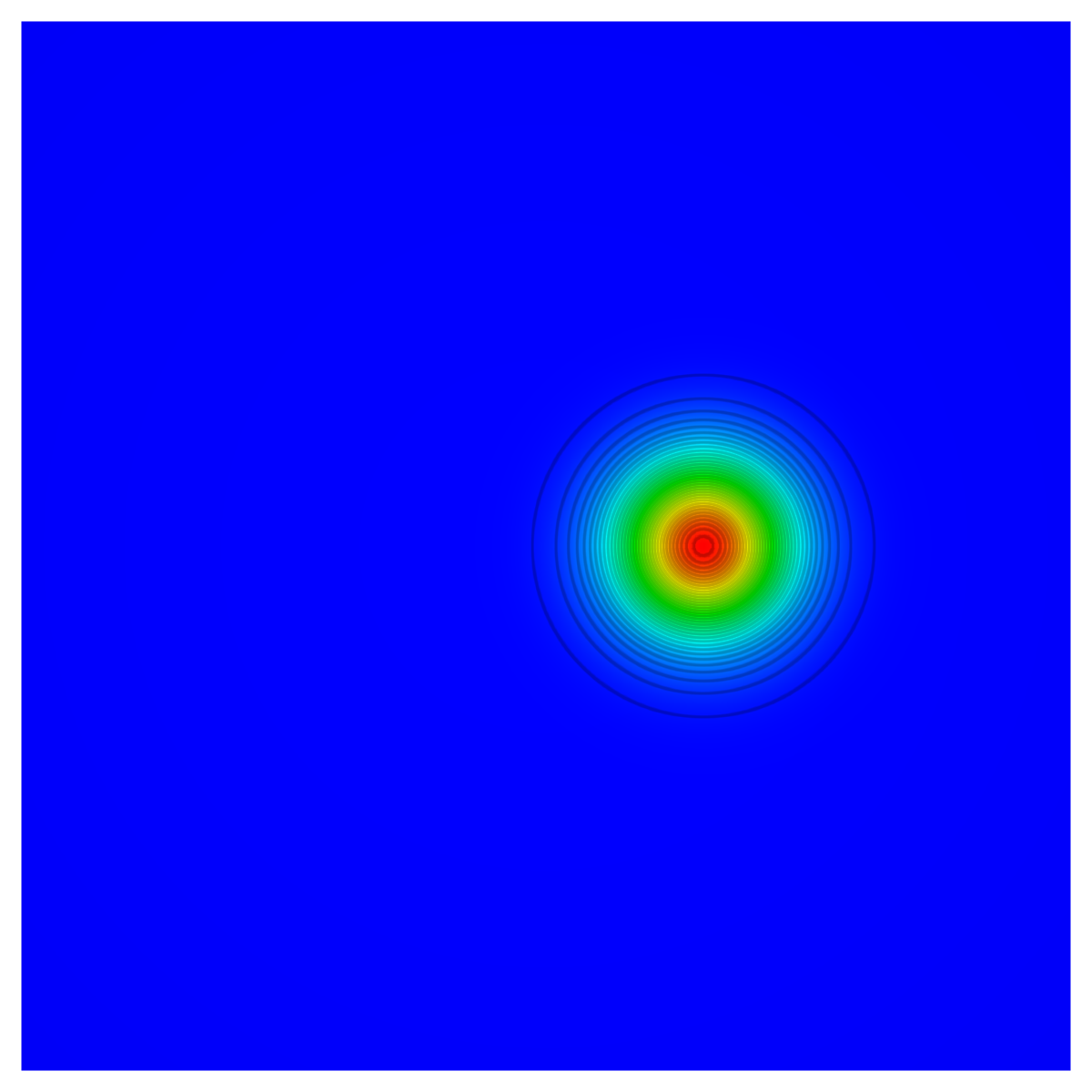}\newline
    $t = 15$\qquad\qquad $t=16$
    \caption{(Swirling deformation flow) Density at times $t=15$ and $t=16$, obtained with the Bell nodal scheme with $\Delta t=1/16$ and $N=256$ ($h\approx 2.5\times 10^{-2}$).} 
    \label{fig:sdf_rho}
\end{figure}

\begin{figure}
    \centering
    \includegraphics[width=0.45\textwidth]{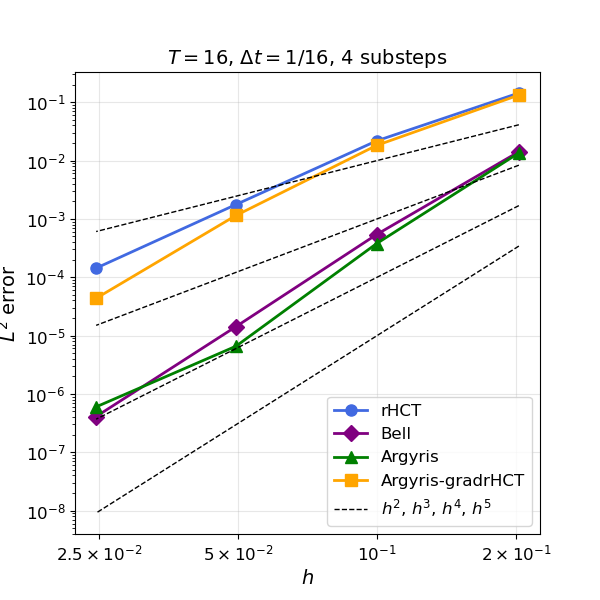}
    \includegraphics[width=0.45\textwidth]{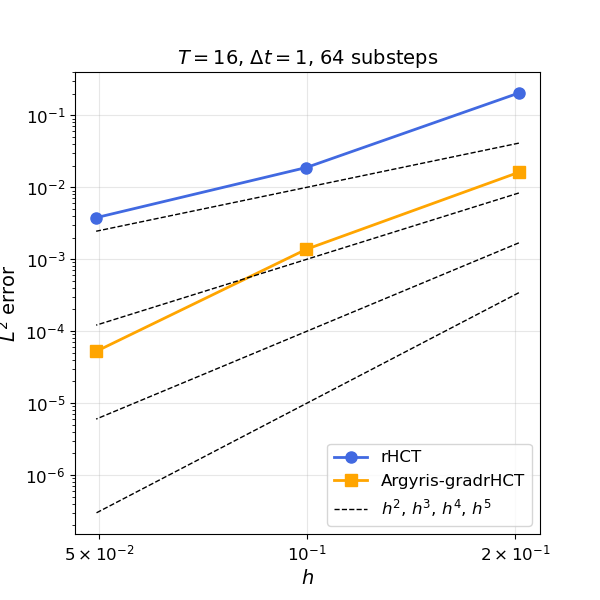}
    \caption{(Swirling deformation flow) Convergence curves with two different time steps. Left: $\Delta t = 1/16$. Right: $\Delta t = 1$.} 
    \label{fig:sdf_convergence}
\end{figure}

\begin{figure}
    \centering
    \includegraphics[width=0.45\textwidth]{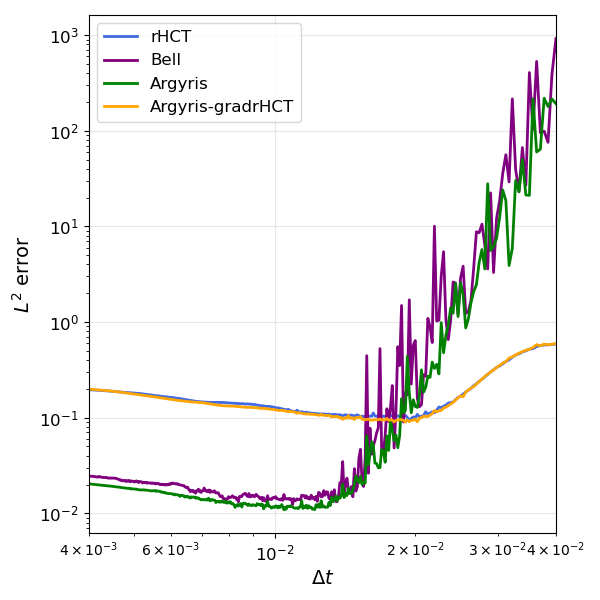}
    \caption{(Swirling deformation flow) $L^2$ error at final time $T= 16$ as a function of the \reviewerA{time step} for nodal SL schemes with $N=32$ ($h \approx 0.2$).} 
    \label{fig:stability_sdf}
\end{figure}

We finally consider the swirling deformation flow test case, whose advection field is defined by:
\begin{equation*}
    a(t,x) = \begin{bmatrix} \sin(\pi x)^2 \sin(2\pi y)\cos(\pi t/T_{\text{sdf}})\\
        -\sin(\pi y)^2 \sin(2\pi x)\cos(\pi t/T_{\text{sdf}})
    \end{bmatrix}.
\end{equation*}
with \reviewerA{$T_{\text{sdf}} = 2$}. With this advection field, the solution to the advection equation equals the initial condition at time $\reviewerA{T_{\text{sdf}}}$:
$
\rho(\reviewerA{T_{\text{sdf}}}, x) = \rho_{\text{init}}(x)$, and this initial condition is the same as in the previous test case (see Fig.~\ref{fig:sdf_rho} for illustration). We note that the characteristic curves are not known analytically and approximation with a RK4 scheme is used with substepping. Figure~\ref{fig:sdf_convergence} shows the convergence curves for two different values of the time step. When considering a small time step, the three considered scheme converge. However, for a larger time step ($\Delta t = 1$), \reviewerA{corresponding to the time interval during which the flow-induced deformation is greatest,} the Bell and Argyris schemes do not work anymore: this is due to the large value of the derivatives of the characteristic field. This is also visible in Figure~\ref{fig:stability_sdf} where \reviewerA{large time steps (small number of iterations)} leads to instabilities for these two schemes. Instead, Argyris-gradrHCT provides better results than rHCT for large time steps, with one order of accuracy larger for fine meshes.

\section{Conclusion}

In this work, we studied high-order Hermite semi-Lagrangian methods for the advection equation on a triangular mesh. With the current implementation, projection SL schemes are much more computationally demanding than the nodal SL schemes. The numerical results show that the 6th order accurate Argyris scheme achieves better accuracy. However the scheme does not work when considering large time steps. In that case, the mixed strategy Argyris-grad rHCT could be advantageous compared with rHCT: although more computationally demanding, it provides slightly better results as being one order more accurate. The future directions of research will be the extension of the methods to non-linear advection equations and the coupling with mesh adaptation.

\paragraph{Acknowkledgments}
Centre de Calcul Intensif d’Aix-Marseille is acknowledged for granting access to its high performance computing resources. This work was conducted within the framework of the ANR project COOKIE, to which the second author
contributed through Task I.3. H. Guillard and B. Nkonga are also acknowledged for invitation 
and fruitful discussions at Saint Etienne de Tinée

%

%
%
\bibliographystyle{abbrv}
\bibliography{biblio}

\section*{Appendix: implementation details}
\addcontentsline{toc}{section}{Appendix}

\paragraph{Projection SL scheme: a Firedrake implementation.} A Firedrake implementation of the scheme is considered, which makes it easy to deal with projection issues. We first describe a classical Discontinuous Galerkin semi-Lagrangian (SLDG) scheme. Starting from a \verb+mesh+, we define a first finite element space \verb+V+, which is a DGd space with $d \in \N$ and project the initial function on it. In Firedrake this writes as follows:
\begin{verbatim}
V = FunctionSpace(mesh, "DG", 3)
f = Function(V).project(f_init)
\end{verbatim}
We then define the mesh \verb+mesh_bwd+, whose vertices are the foot of the characteristic curves originating from the initial mesh nodes, and the associated finite element space \verb+V_bwd+, which is the same DGd space but on this moved mesh. Then each iteration of the SLDG scheme consists into projecting the solution onto the \verb+V_bwd+ and then copy back the degrees of freedom. This reads as follows in Firedrake for the translation equation:
\begin{verbatim}
mesh_bwd = Mesh(mesh.coordinates)
mesh_bwd.coordinates.dat.data_wo[:] -= dt
V_bwd = FunctionSpace(mesh_bwd, "DG", 3)
for it in range(num):
    f_bwd = Function(V_bwd).project(f)
    f.dat.data[:] = f_bwd.dat.data[:]
\end{verbatim}
In order to treat with other spaces, we just have to change the definition of the function space \verb+V+ (and \verb+V_bwd+ accordingly). This is possible for CGd in the current implementation of Firedrake, but not for the Hermite elements (rHCT, HCT, Bell, Argyris). In order to deal with such elements, we first project to the corresponding DG space, then do the remapping and then remap from the DG space to the \verb+V+ space. For HCT and rHCT, we note that this process introduces a modification of the scheme, since the function inside a cell is not polynomial, but piecewise polynomial. These projections also strongly deteriorate the efficiency of the corresponding implementation.

\paragraph{Derivatives of the discrete characteristic curves.} As required in the nodal SL scheme, the derivatives of the characteristic curves are given by:
\begin{align*}
    &\partial_{k}\vc{X}_m(t^n;t^{n+1}, \vc{x})=  \partial_k \vc{x}_m  - \Delta t \sum_{i=1}^p b_i \partial_k \vc{k}_{i,m}(\vc{x}), \quad  \forall k \in \{1, 2\}\\
    &\partial_k\vc{k}_{i,m}(\vc{x}) = \sum_{\ell = 1}^d \partial_\ell \vc{a}_m \Big(\vc{y}_i(\vc{x})\Big)\ \Big(\partial_k \vc{x}_{\ell} - \Delta t \sum_{j=1}^{i-1} a_{i,j} \partial_k \vc{k}_{j,\ell}(\vc{x})\Big),\quad  \forall i \in \interv{1, p},\\
    \text{with}\qquad &\vc{y}_i(\vc{x}) = \vc{x} - \Delta t \sum_{j=1}^{i-1} a_{i,j} \vc{k}_j(\vc{x}),
\end{align*}
where the dependence of the velocity field on time has been omitted for concision. The second partial derivatives are given by:
\begin{align*}
    &\partial_{k'k}\vc{X}_m(t^n;t^{n+1}, \vc{x}) =  - \Delta t \sum_{i=1}^p b_i \partial_{k'k} \vc{k}_{i,m}(\vc{x}),\hspace{3.5cm}\forall k,k' \in \{1, 2\},\\
    &\partial_{k'k}\vc{k}_{i,m}(\vc{x}) = \hspace{8cm} \forall i \in \interv{1, p},\\
    &\sum_{\ell = 1}^d  \sum_{\ell' = 1}^d \partial_{\ell'\ell} \vc{a}_m\Big(\vc{y}_i(\vc{x})\Big) \left(\partial_{k'} \vc{x}^{\ell'} - \Delta t \sum_{j=1}^{i-1} a_{i,j} \partial_{k'} \vc{k}_j^{\ell'}(\vc{x})\right)
    \left(\partial_k \vc{x}_{\ell} - \Delta t \sum_{j=1}^{i-1} a_{i,j} \partial_k \vc{k}_{j,\ell}(\vc{x})\right)\\
    &+ \sum_{\ell = 1}^d \partial_\ell \vc{a}_m \Big(\vc{y}_i(\vc{x})\Big)\left(\Delta t \sum_{j=1}^{i-1} a_{i,j} \partial_{k'k} \vc{k}_{j,\ell}(\vc{x})\right).
\end{align*}

\end{document}